\documentclass[smallextended,12pt]{svjour3}
\smartqed 
\usepackage{amssymb,amsmath,multirow,array,graphicx}
\usepackage{natbib}
\usepackage{breqn}
\usepackage[left=1.5cm, right=1.5cm, top=2.5cm]{geometry}

\begin{document}
%\begin{frontmatter}
\title
{ Reduction of kinetic equations to Li\'enard-Levinson-Smith form: Counting limit cycles}
\author{Sandip Saha$^\dagger$         \and
        Gautam Gangopadhyay$^{\dagger *}$ \and %etc.
        Deb Shankar Ray$^\ddagger$
}

\institute{$^\dagger$S N Bose National Centre For Basic Sciences \at
              Block-JD, Sector-III, Salt Lake, Kolkata 700106, India \\
              \email{gautam@bose.res.in (Corresponding author$^*$)} 
           \and
           $^\ddagger$Indian Association for the Cultivation of Science \at
              Jadavpur, Kolkata 700032, India
}
\date{\today}
\maketitle
\begin{abstract}
We have presented an unified scheme to express a class of  system of equations in two  variables into a   Li\'enard - Levinson - Smith(LLS) oscillator form. We have derived the condition for limit cycle with special reference to Rayleigh and Li\'enard systems for arbitrary polynomial functions of damping and restoring force. Krylov-Boguliubov(K-B) method is implemented to determine the  maximum number of limit cycles admissible for a LLS oscillator atleast in the weak damping limit. Scheme is illustrated by a number of model systems with single cycle as well as the multiple cycle cases.
\end{abstract}

%\keywords{Limit Cycle, Li\'enard and LLS system, Krylov-Boguliubov(K-B) method, Chemical oscillation, Isochronous Oscillator}

\maketitle
\section{Introduction}
Various open kinetic systems\cite{murraybio,murraynld,epstein,strogatz,goldbook,
remickens} in physics, chemistry and biology, are generically described by a minimal model of autonomous coupled differential 
equations\cite{slross,arnold,birkhoff,jjstoker,smale1967,len0} of two variables. They exhibit self sustained oscillation in the form of  stable limit cycle in a phase plane in many examples, such as, chemical reactions\cite{epstein,brusselator2009,len4}, biological rhythms\cite{murraynld,murraybio,kaiser83,kaiser91,k-dsr,strogatz}, vibrations in mechanical\cite{rand2012}, optical system and musical instrument\cite{strogatz,dsrrayleigh}, to name a few. A Rayleigh\cite{dsrrayleigh} equation in violin string and van der Pol oscillation in electric circuit are the classic examples, in this context. More generally,  Li\'enard\cite{perko,len4,strogatz,lls1,lls2,lls3} equation underlines the  concrete criteria for the existence of at least one limit cycle for a general class of such systems of which van der Pol is a special case of the form $\ddot{x}+f(x)\dot{x}+x=0$ where $f(x)=\epsilon (x^2-1)$ and Li\'enard transformation is $\dot{x}=y-F(x)$ and $\dot{y}=-x$ with $F(x)=\int_0^x f(\tau) d\tau$. A further generalisation of Li\'enard equation is the LLS equation\cite{lls1,lls2,lls3,remickens}, $\ddot{x}+F(x,\dot{x}) \dot{x}+G(x)=0$, sometimes called the generalised Li\'enard equation. Casting a general system of kinetic equations  in two variables which describe a variety of scenarios in physical, chemical, bio-chemical and ecological sciences into LLS form\cite{remickens} is often not straight-forward\cite{remickens,len4,limiso,strogatz}. To this end we have provided a scheme for a wide class of open nonlinear equations, cast in the LLS form so that the later becomes amenable to several techniques in nonlinear dynamics. 

Our next objective is to find the nature and the number of limit cycles for a given LLS equation thereby addressing the second part of the Hilbert's $16^{th}$ problem.  The  problem of counting limit cycle has a long legacy since Hilbert, Smale and many others and still continues it without complete understanding\cite{birkhoff,jjstoker,smale1998,perko,gaiko2008,countinglcjkb,infdampinglcjkb,lcbounestjkb}. Our scheme is based on the KB method of averaging\cite{kbbook,strogatz,slross,len0}, a variant of multi-scale perturbation technique\cite{strogatz,chen1,chen2,sarkar2012} to derive amplitude equation with considering the polynomial forms of the nonlinear damping and restoring force functions. We have illustrated our results on a variety of known model systems\cite{perko,gaiko2008,giacomini} with single and multiple limit cycles\cite{perko,gaiko2008,giacomini}. 

\section{Reduction of Kinetic Equations to Li\'enard - Levinson - Smith(LLS) form: Conditions for limit cycle}
We  consider here a set of two dimensional  autonomous kinetic equations for an open system. Our aim is to cast the equations into a form of a variant of LLS oscillator\cite{strogatz,len4,limiso} or LLS oscillator\cite{strogatz,len4,limiso,len0,remickens} which can further be reduced to Rayleigh and Li\'enard form.  Let us begin with the system of autonomous kinetic equations
\begin{align}
\frac{dx}{dt} &= a_0+a_1 x+a_2 y+f(x,y), \nonumber\\
\frac{dy}{dt} &= b_0+b_1 x+b_2 y+g(x,y),
\label{eq1}
\end{align}
where $x(t)$ and $y(t)$ are, for example, field variables or populations of species of chemical, biological or ecological process \cite{strogatz,goldbook,murraynld,murraybio}  with $a_i,b_i$ for $i=0,1,2$ are all real parameters expressed in terms of the appropriate kinetic constants. Let, ($x_s, y_s$) be the fixed point of the system and $f(x,y)$ and $g(x,y)$ are the non-linear functions of $x$ and $y$. A first step is shifting the steady state ($x_s, y_s$) to the origin ($0, 0$) with the help of a linear transformation as LLS system is a second order homogeneous ordinary differential equation. 
%We have observed that the function $f(x,y)$ must have a non-fractional form in $x$ and $y$ to have a Li\'enard or LLS form through linear transformation and there must exist a proportionality relationship between $f(x,y)$ and $g(x,y)$ i.e.  $g(x,y) \propto f(x,y) \implies g(x,y) =\mu f(x,y), \mu\in\mathbb{R}$. This is an important relation to have the Li\'enard form otherwise the linear transformation is not applicable and if we go through non-linear transformation then there may appear some singularities. %%%%%%%%%%%%%%%%%%%
%For a Li\'enard or LLS oscillator, $(0,0)$ is always a fixed point as it is the form of a second order homogeneous ordinary differential equation. The usual approach in taking perturbation around the fixed point\cite{len4,limiso} when casting the kinetic equations into LLS form is basically the method of the conversion of the non-zero fixed point to a zero fixed point whereby stability analysis near zero fixed point becomes trivial.%For systems having $(0,0)$ fixed point need not be applicable by this method provided there must exist a linear transformation like, $\xi=x$, $ \dot{\xi}=\dot{x}=y$ as they can be cast into Li\'enard form through one step by taking derivative, as for example, van der Pol\cite{limiso} system. 

The linear transformation can be chosen by introducing a new pair of variables $(\xi,u)$, both of which are functions of  $x$ and $y$ where $\xi=\beta_0+\beta_1 x+\beta_2 y$ with $\beta_0=-(\beta_1 x_s+\beta_2 y_s)$ i.e. $\xi=\beta_1 (x-x_s)+\beta_2 (y-y_s)$ such that $\dot{\xi}=u$. $\beta_1,\beta_2$ are weighted constants such that it makes the new steady state at the origin, $\xi_s=0$, $u_s=0$. $u$ is expressed as $u=\alpha_0+\alpha_1 x+\alpha_2 y$, with  $\beta_i$, $\alpha_i$ for $i=0,1,2$ are all real constants which can be expressed in terms of system parameters. From the inverse transformation we can easily obtain the expressions for $x$ and $y$ as given by
\begin{align}
x &=\frac{\alpha _2 \left(\beta _0-\xi \right)+\beta _2 \left(u-\alpha _0\right)}{\alpha _1 \beta _2-\alpha _2 \beta _1}=L(\xi,u), \nonumber \\
y &=\frac{\alpha _1 \left(\xi -\beta _0\right)+\beta _1 \left(\alpha _0-u\right)}{\alpha _1 \beta _2-\alpha _2 \beta _1}=K(\xi,u),
\label{eq2}
\end{align}
provided that $\alpha_1 \beta_2-\alpha_2 \beta_1 \neq 0$. Differentiating again, $\dot{\xi}=u$ with respect to the independent variable $t$  we get,
\begin{align}
\ddot{\xi}&=\dot{u} =\alpha_1\dot{x}+\alpha_2\dot{y}\nonumber\\
&=\alpha_1 \lbrace a_0+a_1 L(\xi,\dot{\xi})+a_2 K(\xi,\dot{\xi})+\varphi(\xi,\dot{\xi})\rbrace+\alpha_2 \lbrace b_0+b_1 L(\xi,\dot{\xi})+b_2 K(\xi,\dot{\xi})+\phi(\xi,\dot{\xi})\rbrace,
\label{lls}
\end{align}
where, $L(\xi,\dot{\xi})=c_1 \xi+c_2\dot{\xi}+c_{L}$ and $K(\xi,\dot{\xi})=c_3 \xi+c_4\dot{\xi}+c_{K}$ with
\(
\begin{bmatrix}
c_1 & c_2 & c_{L}\\
c_3 & c_4 & c_{K}
\end{bmatrix}
\)
=$\frac{1}{\alpha_1 \beta_2-\alpha_2 \beta_1}$
\(
\begin{bmatrix}
-\alpha_2 & \beta_2 & \alpha_2 \beta_0-\alpha_0 \beta_2\\
\alpha_1 & -\beta_1 & \alpha_0 \beta_1-\alpha_1 \beta_0
\end{bmatrix}
\).
%Note that none of the constants, $c_{L}, c_{K}$ is negligible and $\alpha_1 \beta_2-\alpha_2 \beta_1 \neq 0$. 
The functions $\varphi$ and $\phi$  can be expressed as a power series expansion as, 
\begin{equation}
\varphi(\xi,\dot{\xi})=\sum_{n,m=0}^{\infty} \varphi_{nm} \xi^n \dot{\xi}^m \hspace{0.5cm} and \hspace{1cm}\phi(\xi,\dot{\xi})=\sum_{n,m=0}^{\infty} \phi_{nm} \xi^n \dot{\xi}^m,
\label{eq4}
\end{equation} 
with, $\phi(\xi,\dot{\xi})=\mu \varphi(\xi,\dot{\xi})$, as the functions $f$ and $g$ are related through $\mu$ by $g=\mu f$, $\mu \in \mathbb{R}$. So, after putting the above form in equation \eqref{lls} one can find, 
\begin{align}
\ddot{\xi}&=\alpha_1 a_0+\alpha_1 a_1(c_1 \xi+c_2 \dot{\xi}+c_{L})+\alpha_1 a_2(c_3 \xi+c_4 \dot{\xi}+c_{K})
+(\alpha_1+\mu \alpha_2) \sum_{n,m=0}^{\infty} \varphi_{nm} \xi^n \dot{\xi}^m \nonumber\\
&+\alpha_2 b_0+\alpha_2 b_1(c_1 \xi+c_2 \dot{\xi}+c_{L})+
\alpha_2 b_2(c_3 \xi+c_4 \dot{\xi}+c_{K}), i.e., \nonumber\\
\ddot{\xi}&=A_{00}+\left(A_{10}+\sum_{n>1} A_{n0} \xi^{n-1}\right) \xi+ \left(A_{01}+\sum_{n>0} A_{n1} \xi^n +\sum_{n\ge0} \sum_{m>1} A_{nm} \xi^n \dot{\xi}^{m-1}\right) \dot{\xi},
\label{eq9}
\end{align}
where, $ \alpha_1 a_0+\alpha_2 b_0+(\alpha_1+\mu \alpha_2) \varphi_{00}+(\alpha_1 a_1 +\alpha_2 b_1) c_{L}+(\alpha_1 a_2 +\alpha_2 b_2) c_{K}=A_{00}=0$ (by definition of a zero fixed point of $\xi$), $A_{10}=\alpha_1(a_1 c_1+a_2 c_3)+\alpha_2 (b_1 c_1+b_2 c_3)+(\alpha_1+\mu \alpha_2) \varphi_{10}$, $A_{01}=\alpha_1(a_1 c_2+a_2 c_4)+\alpha_2 (b_1 c_2+b_2 c_4)+(\alpha_1 +\mu \alpha_2) \varphi_{01}$, $A_{n0}=(\alpha_1+\mu \alpha_2) \varphi_{n0}$, $A_{n1}=(\alpha_1+\mu \alpha_2) \varphi_{n1}$ and $A_{nm}=(\alpha_1+\mu \alpha_2) \varphi_{nm}$, where indices follow the values as given in the summation over $m,n \in \mathbb{Z}^+$. Finally, the above equation looks like, 
\begin{equation}
\ddot{\xi}+F(\xi,\dot{\xi}) \dot{\xi} +G(\xi)=0,
\label{eq12}
\end{equation}
where, the functions $F(\xi,\dot{\xi})$ and $G(\xi)$ are given by
\begin{align}
F(\xi,\dot{\xi})&=-[A_{01}+\sum_{n>0} A_{n1} \xi^n +\sum_{n\ge0} \sum_{m>1} A_{nm} \xi^n \dot{\xi}^{m-1}],   \nonumber\\
G(\xi) &=-[A_{10}+\sum_{n>1} A_{n0} \xi^{n-1} ] \xi.
\label{eq99}
\end{align}
Equation \eqref{eq12} is a well known equation of generalised Li\'enard form called LLS equation. The condition for existence of having at least a locally stable limit cycle of the dynamical system is $F(0,0)<0 \implies A_{01}>0$. It can be shown from the linear stability analysis that there is a relation between $F(0,0)$ and eigenvalues ($\lambda_{\pm}$) with, $F(0,0)=-2\ Re(\lambda_{\pm})$. %For a LLS system, there are some additional conditions\cite{remickens} to have a limit cycle as, (i) $\xi G(\xi)>0$ for $\lvert \xi \rvert\ >0$, (ii) $\int_{0}^{+\infty} G(\xi) d\xi=\int_{0}^{-\infty} G(\xi) d\xi=\infty$, (iii) there exists some $\xi_0 >0$ such that $F(\xi,v) \ge 0$ for $\lvert \xi \rvert \ge \xi_0$, (iv) there exists an $A$ such that for $\lvert \xi \rvert \le \xi_0$, $F(\xi,v) \ge-A$ and (v) there exists some $\lvert \xi_1 \rvert > \xi_0$ such that $\int_{\xi_0}^{\xi_1} F(\xi,v) d\xi \ge 10\ A\ \xi_0$, where $v>0$ is an arbitrary decreasing positive function of $\xi$, which are well established by LLS, explained in \cite{remickens}. 
For a LLS system, there are six conditions to have a limit cycle are given in \cite{lls2,lls3,remickens,len4}.
Out of these six conditions, the condition $F(0,0)$ plays an important role to have a locally stable or unstable limit cycle for such kind of system\cite{len4,limiso,len0,remickens} depending upon the sign of $F(0,0)$ is $<0$ or $>0$, respectively. In particular, two situations may arise:

\textbf{I:} For $A_{nm}=0$, with $n\ge2, \forall m$ i.e. there be an unique steady state ($\xi_s=0$) with restoring force linear in $\xi$, then the above form of \eqref{eq12} looks like 
\begin{equation}
\ddot{\xi}+F_R(\dot{\xi}) \dot{\xi} +G_R(\dot{\xi})\xi=0,
\label{}
\end{equation}
where, 
\begin{align}
F_R(\dot{\xi})&=-[A_{01}+\sum_{m>1} A_{0m} \dot{\xi}^{m-1}],\hspace{10 pt} G_R(\dot{\xi}) =-[A_{10}+\sum_{m>0} A_{1m} \dot{\xi}^{m}],
\end{align}
which is in the form of generalised Rayleigh oscillator\cite{dsrrayleigh}, % with  considering some further approximations analysed in \cite{dsrrayleigh}, 
the limit cycle condition modifies to, $F_R(0)<0$. %Note that, the quotient restriction of $x,y$ for $f$ is not necessary in this case.

\textbf{II:} For $A_{nm}=0$, with $m\ge2,\forall n$, which corresponds to Li\'enard equation with an unique steady state ($\xi_s=0$). This is of the form 
\begin{equation}
\ddot{\xi}+F_L(\xi) \dot{\xi} +G_L(\xi)=0,
\label{}
\end{equation}
where, 
\begin{align}
F_L(\xi)& =-[A_{01}+\sum_{n>0} A_{n1} \xi^n ],\hspace{10pt} G_L(\xi)=-[A_{10}+\sum_{n>1} A_{n0} \xi^{n-1} ] \xi,
\end{align}
where the limit cycle condition is $F_L(0)<0$. We know that, for a Li\'enard system, the damping force function,  $F_L(\xi)$ and the restoring force function, $G_L(\xi)$ are even and odd functions of $\xi$, respectively.

However, for generalised Li\'enard or LLS system the odd-even properties of $G(\xi)$ and $F(\xi,\dot{\xi})$ have complex ramifications\cite{remickens} for practical systems. Here, we have examined the properties with the help of  Krylov-Boguliubov averaging method. 

%We have observed that the function $f(x,y)$ must have a non-fractional form in $x$ and $y$ to have a Li\'enard or LLS form through linear transformation and there must exist a proportionality relationship between $f(x,y)$ and $g(x,y)$ i.e.  $g(x,y) \propto f(x,y) \implies g(x,y) =\mu f(x,y), \mu\in\mathbb{R}$. This is an important relation to have the Li\'enard form otherwise the linear transformation is not applicable and if we go through non-linear transformation then there may appear some singularities.

%For a Li\'enard or LLS oscillator, $(0,0)$ is always a fixed point as it is the form of a second order homogeneous ordinary differential equation. The usual approach in taking perturbation around the fixed point\cite{len4,limiso} when casting the kinetic equations into LLS form is basically the method of the conversion of the non-zero fixed point to a zero fixed State Bank of India - Kalyani Branchpoint. We have replaced it within the very first step by tuning the linear transformation, which gives the origin is always an unique critical point.

\section{Maximum Number of Limit Cycles}
\label{max-cycle}
We now restrict ourselves to the case of LLS systems where the $F(\xi,\dot{\xi})$ and $G(\xi)$ are the polynomial functions of $\xi$ and $\dot{\xi}$. It is well known that linear functional forms of $F$ and $G$ preclude the existence of limit cycle. This can be readily seen by considering the typical examples, e.g., a Harmonic oscillator or a weakly nonlinear oscillator with a potential $\frac{1}{2} \omega_0^2 x^2+\frac{1}{3} \lambda x^4$, $0 < \lambda <1$ or a Lotka-Volterra model, where one encounters a center. We therefore consider the polynomial form of nonlinear damping function $F(\xi,\dot{\xi})$ and restoring force function $G(\xi)$ for our analysis of limit cycle. In what follows we employ K-B method of averaging to show that the characteristic even/odd powers of polynomials play crucial role in determining the behaviour of the associated amplitude and phase equations.

To begin with we consider some fixed values of $m,n$ of equation \eqref{eq99} to truncate the series at $M,N$, for the highest power of $\dot{\xi}$ and $\xi$, respectively. For explicit structure of a prototypical example of an amplitude equation we choose upto $M=N=3$ for  illustration. This includes all possible cases for the even and odd nature of $F(\xi,\dot{\xi})$ and $G(\xi)$, respectively. Then the above form of $F(\xi,\dot{\xi})$ and $G(\xi)$ will be in the following reduced forms,
\begin{align}
F(\xi,\dot{\xi}) &=-[A_{01}+A_{11} \xi+A_{21} \xi^2+A_{31} \xi^3 +A_{02} \dot{\xi}+A_{12} \xi \dot{\xi}+A_{22} \xi^2 \dot{\xi}+A_{32} \xi^3 \dot{\xi} \nonumber\\
&+A_{03} \dot{\xi}^{2}+A_{13} \xi \dot{\xi}^{2}+A_{23} \xi^2 \dot{\xi}^{2}+A_{33} \xi^3 \dot{\xi}^{2}], \nonumber\\
G(\xi) &=-[A_{10} \xi+A_{20} \xi^{2}+A_{30} \xi^{3}].
\end{align}
Let us take $|F(0,0)|=\sigma \in \mathbb{R}^+$, an arbitrary constant with $F(\xi,\dot{\xi})=\sigma F_{\sigma} (\xi,\dot{\xi})$. Then the LLS equation can be rewritten as
\begin{align}
\ddot{\xi}+\sigma F_{\sigma} (\xi,\dot{\xi}) \dot{\xi}+G(\xi)=0.
\label{eq3}
\end{align}
Therefore the final equation takes the form of a non-linear oscillator after rescaling  $t$ by $\tau$ taking, $\omega t \rightarrow \tau$ as
\begin{align}
\ddot{Z}(\tau)+\epsilon h(Z(\tau),\dot{Z}(\tau))+Z(\tau)=0,
\label{14}
\end{align}
where, $0<\epsilon=\frac{\sigma}{\omega^2}\ll 1$, 
$\omega^2=-A_{10}>0$ and $Z(\tau) \equiv \xi(t)$ and $\omega \dot{Z}(\tau) \equiv \dot{\xi}(t)$. Equation \eqref{14} is now ready for the treatment using K-B method with 
\begin{dmath}
h(Z,\dot{Z})=-\left[\lbrace B_{01}+ B_{11} Z+B_{21} Z^2+B_{31} Z^3 + B_{02} \omega \dot{Z} + B_{12} Z \omega \dot{Z}+ B_{22} Z^2 \omega \dot{Z} + B_{32} Z^3 \omega \dot{Z}+ B_{03} \omega^2 \dot{Z}^{2}+B_{13} Z \omega^2 \dot{Z}^{2}+ B_{23} Z^2 \omega^2 \dot{Z}^{2}+B_{33} Z^3 \omega^2 \dot{Z}^{2} \rbrace \omega \dot{Z}+ B_{20} Z^{2}+B_{30} Z^{3}\right], 
\end{dmath}
where $B_{ij}=\frac{A_{ij}}{\sigma}$, $i,j=0,1,2,3$ with $B_{00}=0$ and $B_{01}$ will take the fixed value, -1, 0, or 1 depending upon the nature of the fixed point is stable focus, center/center-type or limit cycle, respectively. Now choosing, $Z(\tau) \approx r(\tau) \hspace{0.1 cm} \cos (\tau+\phi(\tau))$ as a solution of eq. \eqref{14} we have $\dot{Z}(\tau) \approx -r(\tau) \hspace{0.1 cm} \sin (\tau+\phi(\tau))$  with slowly varying radius $r(\tau)=\sqrt{Z^2 +\dot{Z}^2}$ and phase $\phi(\tau)=-\tau+ tan^{-1} (- \frac{\dot{Z}}{Z})$. The function $h(Z,\dot{Z})$ contains all the non-linear terms and $\epsilon$ is the  non-linearity controlling parameter i.e. one has to satisfy $0<\sigma\ll\omega^2$. Then one can obtain $\dot{r}(\tau)=\epsilon h \sin (\tau+\phi(\tau))$ and $\dot{\phi}(\tau)=\frac{\epsilon h}{r(\tau)} \cos (\tau+\phi(\tau))$ i.e. the time derivative of  amplitude and phase are of $O(\epsilon)$. So, after taking a running  average\cite{strogatz,remickens,slross} of a time dependent function $U$ defined as, $\overline{U}(\tau) = \frac{1}{2 \pi} \int_{0}^{2\pi} U(s) ds$, one finds, $\dot{\overline{r}} = \langle \epsilon h \sin (\tau+\phi(\tau)) \rangle_\tau$ and $ \dot{\overline{\phi}}= \langle \frac{\epsilon  h}{r(\tau)} \cos (\tau+\phi(\tau))\rangle_\tau$, which gives,
\begin{align}
\dot{\overline{r}} &= \frac{\epsilon  \omega \overline{r} }{16} \lbrace\overline{r}^2 \left(B_{23} \overline{r}^2 \omega ^2+6 B_{03} \omega ^2+2 B_{21}\right)+8 B_{01}\rbrace+O(\epsilon^2), \nonumber\\
\dot{\overline{\phi}} &= -\frac{\epsilon \overline{r}^2}{16} \left(B_{32} \overline{r}^2 \omega ^2+2 B_{12} \omega ^2+6 B_{30}\right)+O(\epsilon^2).
\label{eq10}
\end{align}

Now from a close look  at the equation for $\dot{\overline{r}}$, it is apparent that only even elements of $F(\xi,\dot{\xi})$  appears but none of any elements of $G(\xi)$ is present due to the zero averages of $\sin^{\mu} \cos^{\nu}$ terms with $\mu=1$ and $\nu \in \mathbb{Z}$. 
%The zero averages will come (i) for the Sine function with ($\mu+\nu$) is even or odd i.e. $\forall$ $\nu \in \mathbb{Z}$ with $\mu=2\eta+1,\forall \eta \in \mathbb{Z}$, together with (ii) odd ($\mu+\nu$), $\forall \mu,\nu \in \mathbb{Z}$. 
The non-zero averages arise  only when $\mu,\nu$ both are even i.e. $\mu=2 \eta_1, \nu=2  \eta_2; \eta_1,\eta_2 \in \mathbb{Z}$. Thus, the effect in $\dot{\overline{r}}$ appears only through the even coefficients of $F(\xi,\dot{\xi})$ i.e. by examining the respective variables in the $\dot{\overline{r}}$ equation, we find that only some even coefficients appear for the first order correction. % and one can say that the function $F(\xi,\dot{\xi})$ is even. 
On the other hand $\dot{\overline{\phi}}$ contains only even coefficients of $F(\xi,\dot{\xi})$ which are not in amplitude equation along with odd coefficients of $G(\xi)$ which shows that only odd $G(\xi)$ plays a role here.%$G(\xi)$ must be an odd function.

So, from the equation of $\dot{\overline{r}}$, one finds that there exist at most $4$ non-zero values of $\overline{r}$. If out of the four roots every pair appears as conjugate then  there are three possibilities. The cases are, (i) two different sets of complex conjugate roots giving an asymptotically stable solution, (ii) one pair of complex conjugate roots and two real roots of equal magnitude with opposite sign implying a limit cycle solution having only one cycle and (iii) either four real roots of equal magnitude with opposite sign having double multiplicity  gives a limit cycle solution with only one cycle or two different sets of real roots of equal magnitude with opposite sign, may give limit cycle solution with two different cycles of different radius. The unique zero values of the roots of $\overline{r}$ gives a center or center-type\cite{len3.5} situation. So, in short, the existence of a nonzero real root will provide the radius of the cycle which will be stable or unstable depending on the $-ve$ or $+ve$ sign of $\frac{d \dot{\overline{r}}}{d \overline{r}}$, at $\overline{r}=\overline{r}_{ss}$ and at $\overline{r}_{ss}=0$ $\frac{d \dot{\overline{r}}}{d\overline{r}}>0$ or $<0$ gives the nature of the fixed point.

As an example, for Kaiser model\cite{kaiser83,kaiser91,k-y2007a,k-y2007b,k-y2007c,k-y2010,k-dsr,k-y2012}, there exist three limit cycles for a certain range of $\alpha, \beta$. So, if we choose the parameters, $\alpha$ and $\beta$ from the three limit cycle zone then there exist six real roots with three different pairs i.e., three different radii exist according to three cycles. But, slightly away from the three limit cycle zone, there will exist only a pair of real roots with the same magnitude and other four will appear as a complex conjugate pairs and together produces only a stable limit cycle.

Note that, to have a stable limit cycle solution, one condition must be satisfied i.e. $F(0,0)<0$. But, it fails to give how many limit cycles the system can admit. According to the root finding algorithm one can guess the maximum number of cycles of a LLS system.  The condition $F(0,0)<0$   plays an important role as a check for the existence of atleast one stable limit cycle. But for  2-cycle situations one can have at first the locally unstable limit cycle before locating the outer stable limit cycle and in this situation $F(0,0)>0$.

Based on these considerations we have prepared a table(Table-I) illustrating the possible cases for the maximum number of non-zero real roots or the limit cycles.
\begin{center}
\textbf{Table-I:} Maxumum number of limit cycles for LLS system
\begin{tabular}{ |c|c|c|c|c| } 
\hline
\hline
$N$ & $M$ & $N+M$ &Max. No. of Non-zero & Max. No. of \\
 & &  &Real Roots (Even) & Cycle(s)\\
\hline
%\multirow{3}{4em}
\hline
&&&&\\
Even & Even & Even & $N+M-2=(N-1)+(M-1)$ & $\frac{N+M}{2}-1$\\ 
&&&&\\
Even & Odd & Odd & $N+M-1=(N)+(M-1)$ &$\frac{N+M-1}{2}$\\
&&&&\\ 
Odd & Even & Odd & $N+M-3=(N-2)+(M-1)$ &$\frac{N+M-3}{2}$\\ 
&&&&\\
Odd & Odd & Even & $N+M-2=(N-1)+(M-1)$ &$\frac{N+M}{2}-1$\\
\hline
\end{tabular}
\end{center}

Now if we denote the non-zero real values of $\overline{r}$ as an existence of limit cycles as $\overline{r}$ gives the radius of the cycle where at the same time a pair of conjugate (one $+ve$ and one $-ve$) roots of equal magnitude exists for such kind of LLS systems then out of these two  roots, radius will be measured by the magnitude and each distinct magnitude counts the number of cycles. For example, if there exists six roots, say, $(p,-p)$ occurring twice and $(q,-q)$ occurring once then the number of cycles will be $2$ of radius $p$ and $q$, respectively. So, if there are all real roots occurring once, then the number of cycles will be atmost $\frac{N+M-2}{2}$ or $\frac{N+M-1}{2}$ or $\frac{N+M-3}{2}$. For LLS equation with $N,M$ are the maximum power of $\xi$ and $\dot{\xi}$ respectively, we have performed the K-B analysis numerically for $N=10,M=10$. The result is given in table-II. For Rayleigh system with $N=1$, for all $M \ge 1$, maximum number of limit cycle will be $\frac{M-1}{2}$ or $\frac{M-2}{2}$ for odd or even $M$, respectively. For Li\'enard system with $M=1$ for all $N \ge 1$, the maximum number of limit cycle becomes $\frac{N-1}{2}$ or $\frac{N}{2}$ for odd or even $N$, respectively. The above table is valid for an arbitrary finite polynomials of $F$ and $G$. For the case of arbitrary infinite polynomial\cite{countinglcjkb,infdampinglcjkb,lcbounestjkb} cases maximum number of limit cycles can be stated for finite truncation.

\section{Applications to Some Model Systems}
Here we have examined three classes of physical models where the analysis of the maximum number of limit cycles holds. This connection with the general model system is discussed with polynomial damping and restoring force function.

\subsection{One-Cycle Cases: van der Pol Oscillator, Simple Glycolytic Oscillator, Modified Brusselator Model}

Considering the van der Pol oscillator\cite{strogatz,limiso,slross,kaiser83,epstein,murraynld,powerlaw} with equation, $\ddot{x}+\epsilon (x^2-1) \dot{x}+x=0$ having the weak nonlinearity for $ 0 < \epsilon \ll 1$ produces a locally stable limit cycle with $F(0,0)<0$.  So, if we compare with the general table we have  $N=2$ and $M=1$. This gives a condition for a unique stable limit cycle. 

Next considering the Li\'enard form\cite{limiso,len3.5} of simple Glycolytic oscillator\cite{strogatz,limiso,len3.5,epstein,murraynld,goldbook} as, $$ \ddot{\xi}+\left[(1+a+3 b^2)-2 b \xi-2bk-3 b \dot{\xi} + \xi\dot{\xi} + k\dot{\xi} + \dot{\xi}^2 \right] \dot{\xi}+(a+b^2)\xi=0; a,b>0, k=b+\frac{b}{a+b^2}$$, has a unique stable limit cycle with $F(0,0)<0$\cite{limiso} having $N=1$ and $M=3$. This gives one limit cycle.

Furthermore, considering the Modified Brusselator model having the Li\'enard form\cite{len4,limiso,goldbook}, $$\ddot{\xi}+\left[ -\frac{2 a_1 \xi}{\alpha}-b+\frac{a_1^2}{\alpha^2}+\alpha-\frac{2 a_1 \dot{\xi}}{\alpha^2}+\frac{b \dot{\xi}}{a_1}+\frac{\dot{\xi}^2}{\alpha^2}+\frac{\xi \dot{\xi}}{\alpha} \right] \dot{\xi} +\frac{a_1^2 \xi}{\alpha}=0; a_1,b,\alpha > 0$$ gives a unique stable limit cycle with $F(0,0)<0$\cite{len4,limiso}, where $N=1$ and $M=3$ again giving rise to the same situation.

\subsection{Two-Cycle Cases}

We rewrite the Li\'enard form according to Ref.\cite{perko,giacomini} as, $\dot{x(t)}=y(t)-F(x(t)), \dot{y(t)}=-x(t)$, where $F(x(t))$ is an odd polynomial. After taking derivative it takes the following form $\ddot{x}+F'(x)\dot{x}+x=0$, where $F'(x)=\frac{\partial F(x)}{\partial x}$ now becomes the form of an even polynomial. For $F(x)=a_1 x+a_2 x^2+a_3 x^3$, it has been shown\cite{perko,giacomini,lins} that the system allows a unique limit cycle if $a_1 a_3 <0$, which will be stable if $a_1 < 0$ and unstable if $a_1 > 0$. This corresponds to the table $N=2,M=1$. Further extension by Rychkov\cite{rychkov} shows that for $F(x)=(a_1 x + a_3 x^3 + a_5 x^5)$ the number of limit cycle is atmost two. Numerical simulation corroborates this observation when $F(x)$ is chosen as in Ref.\cite{perko,giacomini,lins},  $F(x)=0.32 x^5 - \frac{4}{3} x^3 + 0.8 x$. For this case the inner one is unstable limit cycle as $F(0,0)=0.8>0$ i.e. a stable fixed point but the outer one is a stable limit cycle. Here, as per Table-I we have $N=4$ and $M=1$. The above table thus gives the strategies to find out the number of limit cycles(both stable and unstable) a system can have. On the other hand our analysis by K-B method provides a hint towards a choice of the parameter space for search of real roots of the radial equation.

\subsection{Three-Cycles Case: Kaiser Bi-rythmicity Model}
Extending  van der Pol oscillator model with a nonlinear function of higher order polynomial,  Kaiser\cite{kaiser83,kaiser91,k-y2007a,k-y2007b,k-y2007c,k-y2010,k-dsr,k-y2012} has described bi-rythmicity with the nonlinear equation, 
\begin{align}
\ddot{x}-\mu (1-x^2+\alpha x^4-\beta x^6) \dot{x}+x=E \cos \Omega t. \label{kaisereq}
\end{align}
Here, $\alpha, \beta, \mu > 0$ tune the non linearity. This is a prototype self-sustained oscillatory system in absence of $E$ and $\Omega$ which are  the amplitude and the frequency of the external excitation, respectively.  The  model exhibits an extremely rich bifurcation behaviour and the system actually produces bi-rhythmicity. It has been emphasized that in the undriven case, the model is a multi-limit cycle oscillator and has three limit cycles, two of them are stable and between the two stable limit cycles there is an unstable one which divides the basins of attraction of the two stable cycles. In presence of $E$ and $\Omega$, the above system exhibits some interesting features\cite{kaiser83,kaiser91,k-y2007a,k-y2007b,k-y2007c,k-y2010,k-y2012}. From table-I one has $N=6$ and $ M=1$ with even-odd sub cases, while $E=0$. Thus, there may have 6 roots for the radial equation if $\mu>0$ and $\alpha,\beta$(controlling parameters of the radii) are chosen from three limit cycle zone($\alpha=0.144,\beta=0.005$) and finally, the number of distinct values will be 3 which  implies that the system can have atmost three limit cycles (but here it is exactly 3). Further, if we choose $\beta=0$ for the above undriven Kaiser model with $\alpha=0.1$, one can have two limit cycles with radii $\approx 2.35$ and $\approx 3.80$, respectively, of which the smaller one will be stable and the larger one will be unstable.

Note that, if there are odd number of limit cycles, say $l$, then out of the $l-$cycles, $\frac{l+1}{2}$ will be stable limit cycles and the remaining $\frac{l-1}{2}$ will be unstable limit cycles. For example, for the van der Pol oscillator, Glycolytic oscillator, Modified Brusselator model etc. only one limit cycle exists which is stable. For Kaiser model, $l=3$ and one can observe the situations accordingly. So, for odd number of cycles innermost one will be locally stable. %one can see the same situation in Glycolytic bi-rythmicity\cite{goldbook} model, which is difficult to cast to the Li\'enard form. But for even number case  one can't distinguish as above and some other situations may arise.

\subsection{k-Cycle Cases: }

\subsubsection{A Model With N=1 and M=2k+1}
For counting the number of limit cycles  Gaiko\cite{gaiko2008} has shown, for a Li\'enard-type system i.e., LLS equation having the form,
\begin{align}
\ddot{x}- \left(\mu_1+\mu_2 \dot{x}+\mu_3 \dot{x}^2+ \dots +\mu_{2 k} \dot{x}^{2k-1}+\mu_{2 k+1} \dot{x}^{2 k} \right)\dot{x}+x &=0,
\label{gaiko}
\end{align}
can have atmost $k$ limit cycles if and only if, $\mu_1>0$. The result\cite{gaiko2008} correlates with our result. For any value of $k \in \mathbb{Z}$ it fits the odd-odd case of the general table and accordingly,  $M$ and $N$ are $2k+1$ and $1$, respectively, and finally the number of cycles will be atmost $\frac{N+M}{2}-1=k$.

\subsubsection{A Model With N=2k and M=1}

Blows and Lloyd\cite{blows,perko} have stated that ``For the Li\'enard or LLS system $\dot{x}=y-F(x), \dot{y}=-g(x)$ with $g(x)=x$ and $F(x)=a_1 x+a_2 x^2+ \dots + a_{2k+1} x^{2k+1}$ has at most $k$ local limit cycles and there are coefficients with $a_1,a_3,\dots,a_{2k+1}$ altering in sign".  This can be found from the table-I with $N=2k$ and $M=1$ to give the condition of atmost $k$ limit cycles. For  example, taking $k=3$ with $F(x)=-\epsilon(72 x-\frac{392}{3} x^3+\frac{224}{5} x^5 - \frac{128}{35} x^7)$ has exactly three limit cycles for sufficiently small $\epsilon \neq 0$ which are  circles with radii $1,2$ and $3$. The above statement nicely corresponds to the Theorem-6, pp-260\cite{perko}.

\begin{center}
\textbf{Table-II:} Table for highest degree polynomial $N+M (\oplus)$ for LLS system together with maximum number of distinct conjugate roots$(R)$, with $1 \le N, M \le 10$

\hspace*{-1.6 cm}
\begin{tabular}{c|c|c|c|c|c|c|c|c|c|c|c|}
 \rotatebox{0}{\tiny{\textbf{$\oplus, R$}}}& \multicolumn{1}{c}{} & \multicolumn{2}{c}{}& \multicolumn{1}{c}{}&\multicolumn{1}{c}{}& \multicolumn{2}{c}{\textbf{\tiny{M}}}& \multicolumn{2}{c}{}& \multicolumn{2}{c}{}\\\cline{1-12}
 
 & \multicolumn{1}{c}{}& \multicolumn{1}{c}{1} & \multicolumn{1}{c}{2}  & \multicolumn{1}{c}{3}& \multicolumn{1}{c}{4} & \multicolumn{1}{c}{5}  & \multicolumn{1}{c}{6}& \multicolumn{1}{c}{7} & \multicolumn{1}{c}{8}  & \multicolumn{1}{c}{9}& \multicolumn{1}{c}{10} \\\cline{3-12}%

%& 0 & $(2,0)$ & $(3,2)$& $(4,2)$& $(5,4)$  & $(6,4)$ & $(7,6)$ & $(8,6)$ & $(9,8)$& $(10,8)$ \\\cline{3-11}
 
& 1 &$ 2,0 $ & $ 3,0 $ & $ 4,2 $& $ 5,2 $& $ 6,4 $  & $ 7,4 $ & $ 8,6 $ & $ 9,6 $ & $ 10,8 $& $ 11,8 $ \\\cline{3-12}
  
& 2 & $ 3,2 $ & $ 4,2 $ & $ 5,4 $& $ 6,4 $& $ 7,6 $  & $ 8,6 $ & $ 9,8 $ & $ 10,8 $ & $ 11,10 $& $ 12,10 $ \\\cline{3-12}
 
& 3 & $ 4,2 $ & $ 5,2 $ & $ 6,4 $& $ 7,4 $& $ 8,6 $  & $ 9,6 $ & $ 10,8 $ & $ 11,8 $ & $ 12,10 $& $ 13,10 $ \\\cline{3-12}
  
& 4 & $ 5,4 $ & $ 6,4 $ & $ 7,6 $& $ 8,6 $& $ 9,8 $  & $ 10,8 $ & $ 11,10 $ & $ 12,10 $& $ 13,12 $& $ 14,12 $ \\\cline{3-12}
  
\multirow{2}*{\textbf{\tiny{N}}}  & 5 & $ 6,4 $ & $ 7,4 $ & $ 8,6 $& $ 9,6 $& $ 10,8 $  & $ 11,8 $ & $ 12,10 $ & $ 13,10 $& $ 14,12 $& $ 15,12 $ \\\cline{3-12}

& 6 & $ 7,6 $ & $ 8,6 $ & $ 9,8 $& $ 10,8 $& $ 11,10 $  & $ 12,10 $ & $ 13,12 $ & $ 14,12 $ & $ 15,14 $& $ 16,14 $ \\\cline{3-12}
  
& 7 & $ 8,6 $ & $ 9,6 $ & $ 10,8 $& $ 11,8 $& $ 12,10 $  & $ 13,10 $ & $ 14,12 $ & $ 15,12 $ & $ 16,14 $& $ 17,14 $ \\\cline{3-12}
 
& 8 & $ 9,8 $ & $ 10,8 $ & $ 11,10 $& $ 12,10 $& $ 13,12 $  & $ 14,12 $ & $ 15,14 $ & $ 16,14 $ & $ 17,16 $& $ 18,16 $ \\\cline{3-12}
  
& 9 & $ 10,8 $ & $ 11,8 $ & $ 12,10 $& $ 13,10 $& $ 14,12 $  & $ 15,12 $ & $ 16,14 $ & $ 17,14 $ & $ 18,16 $& $ 19,16 $ \\\cline{3-12}
 
& 10 & $ 11,10 $ & $ 12,10 $ & $ 13,12 $& $ 14,12 $& $ 15,14 $  & $ 16,14 $ & $ 17,16 $ & $ 18,16 $ & $ 19,18 $& $ 20,18 $ \\\cline{3-12}
\end{tabular}
\end{center}

Counting the number of limit cycles through RG method in first order will give similar result which was  done by Das et. al.\cite{countinglcjkb,infdampinglcjkb,lcbounestjkb} for some models. We have verified  similar results for (3,3) polynomial cases  for ($F$,$G$) functions using RG method which become increasingly very difficult and almost impossible upto (10,10) case than K-B averaging method as tabulated in this work. It is  very useful to count the number of limit cycles  from the table by just looking at the LLS form. For example, the number of limit cycles  of all models in Ref.\cite{countinglcjkb,infdampinglcjkb,lcbounestjkb} along with the models in our work can  be estimated from our table. The table-I can also be utilized to prepare a model of a desired number of limit cycles in a systematic way.

\section{Conclusions}

 We have presented a scheme to cast a set of a class of coupled nonlinear equations in two variables into a LLS form. By expressing the nonlinear damping and forcing functions as polynomial we have implemented K-B method of averaging to explore the number of admissible limit cycles of the dynamical systems. Our results can be summarised as follows: 

%\item Here we have provided a short-cut scheme to express a coupled non linear system of equations into a LLS form through linear transformation. For this approach one needs a proportionality relationship between the nonlinear coupled functions expressed in power series.

%\item From the above K-B analysis, it can be concluded that the damping force function $F(\xi,\dot{\xi})$ will be even in nature and the restoring force $G(\xi)$ will be odd, although the analysis is for the weakly non-linear case. One can convert a strongly nonlinear  LLS system into a weakly nonlinear by rescaling with respect to the system frequency.
\begin{enumerate}
\item For a LLS system, the number of limit cycles will be atmost $\frac{N+M-2}{2}$ when $N$ and $M$ degree of the polynomials for damping and restoring force both are even or odd. Again, $\frac{N+M-1}{2}$ cycles can be found when $N$ is even and $M$ is odd and finally, $\frac{N+M-3}{2}$ cycles when $N$ is odd and $M$ is even.

\item For a Li\'enard system, in particular, the formula of counting the number of limit cycles follows the same with $M=1$ and $N \in \mathbb{Z^+}$. Also for the generalised Rayleigh situation there occurs a linear restoring force so that $N$ is $1$ and $M \in \mathbb{Z^+}$ .%and follow the same formula. One can verify the examples of van der Pol system as Li\'enard and Rayleigh oscillator for the verification of the above statement.
\item We have validated our general result with the help of a variety of physical systems with one, two, three upto arbitrary k-cycles.

\item This method stated in our work can also be utilized to prepare a model of a desired number of limit cycles in a systematic way.
\end{enumerate}
\noindent
{\bf Acknowledgement}

\noindent
{ Sandip Saha acknowledges RGNF, UGC, India for the partial financial support.} \\\\
\noindent
{\bf Compliance with ethical standards\\
Conflicts of interests}

\noindent
{ The authors declare that there is no conflict of interests regarding the publication of this paper.}

\section{Appendix: Lotka-Volterra System}
To obtain the LLS form of Lotka-Volterra System\cite{epstein,strogatz,goldbook}, let us set $z=\delta x+\beta y$ then $\dot{z}=\alpha \delta x-\beta \gamma y=u$ $\implies$ $x=\frac{\dot{z}+\gamma z}{(\alpha+\gamma)\delta}$ and $y=\frac{-\dot{z}+\alpha z}{(\alpha+\gamma)\beta}$. After taking $t$ derivative upon $\dot{z}$ one can have,
\begin{align*}
\ddot{z} &=(\alpha-\gamma)\dot{z}+\alpha \gamma z +\frac{\dot{z}^2}{\alpha+\gamma}+\frac{\gamma - \alpha}{\alpha+\gamma} z \dot{z}-\frac{\alpha \gamma}{\alpha+\gamma} z^2.
\end{align*}
The fixed point $(0,0)$ gives a saddle solution which is not of any interest in the present context. Choosing the remaining non-zero fixed point for further investigations,  and after taking perturbation $z=\xi+z_s$ around the fixed point $z_s=\alpha+\gamma=\delta x_s+\beta y_s \neq 0$, one can get the LLS form with
%
%\begin{align}
%\ddot{\xi}+F(\xi,\dot{\xi}) \dot{\xi}+G(\xi)=0
%\label{eq2}
%\end{align}
%where, 
$F(\xi,\dot{\xi})=a_1 \xi +a_2 \dot{\xi}$ % is linear and having weak nonlinearity in the LLS equation 
with $a_1=\frac{\alpha-\gamma}{\alpha+\gamma}$ and $a_2= - \frac{1}{\alpha+\gamma}$. It is to be noted that $G(\xi)$ contains nonlinearity with $G(\xi)=\omega^2 \xi+a_3 \xi^2$ where $\omega=\sqrt{\alpha \gamma}=Im(\lambda)$($+ve$ sense) and $a_3=\frac{\alpha \gamma}{\alpha+\gamma}$. After introducing a small parameter $\epsilon_1$ (say) in the constants, $a_i, b_i$  such that $a_i=\epsilon_1 b_i, i=1,2,3$ the above equation reduces to $\ddot{\xi}+\epsilon_1 (b_1 \xi+b_2 \dot{\xi}) \dot{\xi}+\omega^2 \xi+\epsilon_1 b_3 \xi^2=0$. %, and after rescaling near the frequency it reduces to $\ddot{Z}(\tau)+\epsilon h(Z (\tau),\dot{Z}(\tau))+Z (\tau)=0$ with $h=(k_1 Z + k_2 \dot{Z}) \dot{Z}+k_3 Z^2$ where $k_1=\omega b_1$, $k_2=\omega^2 b_2$,$k_3=b_3$, and $\epsilon=\frac{\epsilon_1}{\omega^2}$. % and since we are analysing it in the weak nonlinear case then $\epsilon$ must be set in between 0 and 1 i.e. $0<\epsilon\ll1$ i.e. $0<\epsilon_1\ll\omega^2$. 
%Due to smallness of $\epsilon$, we have $0<\epsilon_1\ll\omega^2=\alpha \gamma \le 1$ $\implies \alpha \le \frac{1}{\gamma}$, since $\alpha,\gamma>0$.

%\bibliographystyle{apa}
\bibliographystyle{ieeetr}
\bibliography{References-ssaha}{}

\end{document}